\documentclass[12pt]{article}
\usepackage{latexsym}
\usepackage{amsfonts}
\def\ZZ{\mathbb Z}
\def\RR{\mathbb R}

\def\CC{\mathbb C}

\def\cpb{\overline{\mathbb C \mathbb P}_2}
\def\sxs{[S^1\times S^3]}

\def\bea{\begin{eqnarray*}}
\def\eea{\end{eqnarray*}}

\newtheorem{main}{Theorem}

\newtheorem{thm}{Theorem}
\newtheorem{prop}[thm]{Proposition}
\newtheorem{cor}[thm]{Corollary}

\newtheorem{defn}[thm]{Definition}

\newenvironment{proof}{\medskip \noindent
{\bf Proof.}}{\hfill \rule{.5em}{1em}
\\}
\newenvironment{xpl}{\mbox{ }\\ {\bf  Example}\mbox{ }}{
\hfill $\diamondsuit$\mbox{}\bigskip}

\begin{document}
\sloppy
\title{Curvature,  Connected Sums, and Seiberg-Witten Theory}
 
\date{Final Version}

\author{Masashi Ishida and Claude LeBrun\thanks{Supported 
in part by  NSF grant DMS-0072591.}  
  }

\maketitle

\begin{abstract} We consider several   differential-topological  invariants 
of compact $4$-manifolds  which   directly arise
 from Riemannian variational problems.
Using recent results of Bauer and Furuta \cite{baufu,bauer2}, we   
compute these invariants in many 
cases that were previously intractable. 
In particular, we  are now able to calculate the 
Yamabe invariant for many connected sums of complex surfaces. \end{abstract}

\section{Introduction}\label{intro}

Given  a smooth compact  $n$-dimensional manifold $M$,
the  {\em canonical metric problem} \cite{bercent} would ask us 
 to geometrize $M$ by somehow endowing it with a preferred 
Riemannian metric. In practice, this usually amounts to  trying to find a
critical point for some natural curvature functional on the space of 
all Riemannian metrics. While few such efforts seem destined to succeed,
they often nevertheless  lead to interesting 
diffeomorphism invariants of  compact manifolds. A prototypical
 example of  such an 
invariant, called the {\em Yamabe invariant} or {\em sigma constant},
historically arose from  an  attempt to construct Einstein metrics 
on  
manifolds of dimension $n > 2$.

Einstein metrics (i.e. metrics of constant Ricci curvature) are 
precisely \cite{bes} the critical points of 
the {\em normalized total scalar curvature}
$${\mathfrak   S} (g)=V_g^{(2-n)/{n}} \int_{M}s_{g}d\mu_{g} ,$$
considered as a functional with domain the space of Riemannian metrics $g$ on
a given smooth compact manifold $M$ of dimension $n > 2$; 
here 
  $s={R^{jk}}_{jk}$ denotes the scalar curvature,  
$d\mu$ is the Riemannian volume measure, 
and $V=\int_{M}d\mu$ is the  total 
volume of $M$ with respect to the relevant metric. 
Unfortunately, one cannot hope to 
construct an Einstein metric on $M$ by minimizing or maximizing 
$\mathfrak S$, as the functional is 
neither bounded above nor bounded below. However, as was first 
pointed out  by H.\ Yamabe \cite{yam},
 the restriction of 
${\mathfrak  S}$ to any {\em conformal class}
$$\gamma = [g] = \{ ug ~|~ u: M \stackrel{C^{\infty}}{\longrightarrow}
{\RR}^{+}\}$$ 
of metrics {\em  is} always 
bounded below,
and Yamabe therefore hoped   to   construct  
 Einstein metrics by instead first minimizing $\mathfrak S$ in 
each conformal class, and then maximizing over 
the set of all conformal classes. While this certainly
fails in general, it nevertheless yields a 
real-valued 
diffeomorphism invariant  \cite{okob,sch}
$${\mathcal Y}(M)=
\sup_{\gamma}\inf_{g\in\gamma}
{\mathfrak  S}(g),$$
 called the {\em Yamabe invariant}  of the 
smooth compact manifold $M$. 
It is not hard to show that 
  ${\mathcal Y}(M) > 0$
if and only if 
 $M$ admits a metric of positive scalar curvature; thus, 
the problem of computing the Yamabe invariant 
may be thought of as a quantitative refinement 
of the question of whether a given manifold admits 
positive-scalar-curvature metrics.

Until recently,  however, the only  available  non-trivial 
exact computations of 
 the Yamabe invariant  concerned the underlying $4$-manifolds of  
   complex algebraic surfaces \cite{lno,lky}, or
 their  connected sums \cite{jp2} with copies of $S^1\times S^3$. 
Our primary  objective in the 
 present article  is to   remedy  this dearth of knowledge, 
by proving  the following result: 

\begin{main} \label{refocus}
Let $X_j$, $j=1, \ldots, 4$, be minimal complex surfaces with 
\begin{eqnarray}b_1(X_j)&=&0, \label{tom} \\
b_+(X_j)&\equiv& 3 \bmod 4 ,\label{dick}  \\
{\textstyle \sum_{j=1}^4}b_+(X_j)&\equiv &4\bmod 8.  \label{harry}
\end{eqnarray}
Let $N$ be any smooth, compact orientable $4$-manifold with
$b_+(N) =0$  which admits a metric of scalar curvature  $\geq 0$. 
 Then, for each $m=1, \ldots, 4$,   
$${\mathcal Y}([\#_{j=1}^mX_j]\#N) =-4\pi \sqrt{2 \sum_{j=1}^m c_1^2 (X_j)} .$$
\end{main}

Here $\#$ denotes the {\em connected sum} operation, which 
is carried out by removing  standard balls, and gluing along the 
resulting boundaries. 
The most important single ingredient in the proof is a 
recent breakthrough in  Seiberg-Witten
theory  due to  Bauer  and Furuta \cite{baufu,bauer2}.

The Yamabe invariant ${\mathcal Y}(M)$ is closely related to 
the diffeomorphism invariant 
 defined \cite{bcg1,lric} by
\begin{equation}\label{is}
{\mathcal I}_{s}(M)  = 
 	\inf_{g}\int_{M}|{s}_{g}|^{n/2}d\mu_{g},
\end{equation}
where the infimum is taken over all Riemannian metrics  on $M$,
and where $s$ and $d\mu$  once again respectively denote the scalar curvature and 
the Riemannian volume measure; cf. Proposition \ref{func} below. 
By definition, ${\mathcal I}_s$  is thus a quantitative obstruction to the
existence of scalar-flat metrics. 
 The power 
of $n/2$ is dictated by the need for  scale invariance: for trivial reasons, 
any other choice would yield the zero invariant.  
But, even with this careful choice, 
it  is   hardly  obvious that the above invariant is non-trivial; 
 and indeed, 
Petean \cite{jp3},  building upon 
the earlier work of Gromov-Lawson \cite{gvln} and Stolz
\cite{stolz},
has proved that  ${\mathcal I}_{s}$ vanishes for every  
simply connected 
$n$-manifold with $n \geq 5$.  
Dimension four, however, turns out to be quite exceptional  for questions 
concerning 
 the scalar curvature, and 
 Seiberg-Witten theory \cite{witten} turns out to predict that 
${\mathcal I}_s$ is non-zero for many simply connected $4$-manifolds. 
Indeed, in proving Theorem \ref{refocus} 
we will simultaneously prove  the following:

\begin{main} \label{focus} 
Let $X_j$, $j=1, \ldots, 4$  and $N$ be as in Theorem \ref{refocus}.
Then, for each $m=1, \ldots, 4$, 
$${\mathcal I}_s ([\#_{j=1}^mX_j]\# N) ={\mathcal I}_s ( \#_{j=1}^mX_j)= 
32\pi^2 \sum_{j=1}^m c_1^2 (X_j) .$$
\end{main}

\begin{xpl}
Let $X$ be the ramified  double  cover of ${\mathbb C \mathbb P}_2$, 
with branch locus 
 a smooth complex curve of degree $8$.
Then  the canonical class of $X$ is the 
pull-back  of the hyperplane class from ${\mathbb C\mathbb P}_2$,
and it therefore follows that  $p_g(X)= 3$ and $c_1^2(X)= 2$. 
In particular, $b_+(X)=7\equiv 3 \bmod 4$,  so we have
$${\mathcal Y}(X \# X ) = 
-4\pi \sqrt{2(2+2)}
=- 8\pi \sqrt{2} $$
and 
$${\mathcal I}_s(X\# X ) =   32\pi^2 (2+2)= 128 \pi^2.$$

By contrast, consider the connected sum $X\# \overline{X}$, where $ \overline{X}$
denotes $X$ with the {\em reverse orientation}. Because $X$ was constructed
 as a branched cover  of ${\mathbb C \mathbb P}_2$,  it 
has a handle decomposition without $1$- or $3$-handles, and 
 the non-spin $4$-manifold  $X\# \overline{X}$ 
 consequently dissolves  \cite[Cor. 6.3.19 \& Prop.  9.1.16]{gost}:
$$X\# \overline{X} \stackrel{\mbox{\it \tiny diff}}{\approx} 44 {\mathbb C \mathbb P}_2\# 44 \cpb . $$
In particular, $X\# \overline{X}$ admits  metrics of
positive scalar curvature, and so, by Proposition \ref{func} below,  
$${\mathcal I}_s(X\# \overline{X} ) = 0.$$
We can also  estimate ${\mathcal Y}(X\# \overline{X} )$
on the basis  of this picture, although 
current techniques do not lend themselves to an exact calculation of   this number. 
Indeed,  ${\mathcal Y} ({\mathbb C \mathbb P}_2) =
{\mathcal Y} (\cpb )  > 0$, so that 
 O.\  Kobayashi's estimate  \cite{okob} of the Yamabe invariant of a  connected sum
 tells us that  
$${\mathcal Y}(S^4) \geq 
{\mathcal Y}(k {\mathbb C \mathbb P}_2\# \ell  \cpb )\geq  {\mathcal Y}({\mathbb C\mathbb P}_2),$$
and we therefore see that 
$${\mathcal Y}(X\# \overline{X} ) \in [12 \pi\sqrt{2} , 8\pi \sqrt{6}].$$
This vividly illustrates the point that 
  $\mathcal Y$ and ${\mathcal I}_s$  are remarkably sensitive to the 
relative orientation of the summands in a connected sum! 
\end{xpl}

Of course, curvature comes in many different flavors, and so it is equally  natural to consider
 invariants such as  
\begin{equation}\label{tight}
{\mathcal I}_{r}(M) = 	\inf_{g}\int_{M}|{r}|^{n/2}_{g}d\mu_{g} ,
\end{equation}
where $r$ denotes the Ricci tensor, and where, once again, 
$n=\dim M$. It is easy to show that  there is a 
tautological inequality 
\begin{equation}\label{taut}
 {\mathcal I}_{r}(M) \geq  n^{-n/4} ~  {\mathcal I}_{s}(M),
\end{equation}
and  that  
equality holds if the Yamabe invariant ${\mathcal Y}(M)$ is both non-positive and 
 realized by an Einstein metric; the failure of equality would thus
give a quantitative obstruction to Yamabe's program for 
finding Einstein metrics.
   It is therefore  interesting to 
observe that the 
 inequality (\ref{taut}) is often strict   in dimension $4$. Indeed, we are  able to
prove: 

\begin{main} 
Let $X_j$, $j=1, \ldots, 4$ be minimal complex surfaces satisfying
(\ref{tom}--\ref{harry}).  Let $N$ be a smooth compact oriented
$4$-manifold which admits an anti-self-dual metric of positive scalar curvature. 
Then, for each choice of $m=1, \ldots, 4$,  
$$ 
{\mathcal I}_r ([\#_{j=1}^mX_j]\#N)
=8\pi^2\left[4m - (2\chi + 3\tau ) (N) + \sum_{j=1}^m c_1^2(X_j)\right]. 
$$
\end{main}
Here, as throughout, $\chi$ and $\tau$ respectively denote the 
signature and Euler characteristic of an oriented $4$-manifold.
 Note that the above hypotheses regarding $N$ imply that $b_+(N)=0$, 
and hence  that $(2\chi + 3\tau ) (N)\leq 4$.
In particular, the Yamabe sup-inf on these manifolds is never 
realized by an Einstein metric if  
 $m > 1$.

\begin{xpl}
Let $X$ once again be the ramified  double  cover of ${\mathbb C \mathbb P}_2$
branched over 
 a smooth octic.
Then 
\begin{eqnarray*}
{\mathcal I}_r(X\# X )
&=&{\mathcal I}_r(X\# X \# S^4)\\
&=&8\pi^2\left[4\cdot 2 - (2\chi + 3\tau ) (S^4) + c_1^2(X)+c_1^2 (X)\right]\\
&=& 8\pi^2\left[4\cdot 2 - 4 + 2+2\right] 
\\&=& 64\pi^2.
\end{eqnarray*}
On the other hand, we have already seen that 
$$
{\mathcal I}_s(X\# X )= 128 \pi^2, 
$$
so that 
$${\mathcal I}_r(X\# X )- \frac{1}{4} {\mathcal I}_s(X\# X ) = 32\pi^2.$$
This gap measures the limiting size of the trace-free part of the Ricci tensor 
as one approaches the Yamabe sup-inf. Roughly speaking, the gap  arises because 
the geometry is trying to break apart into two copies of the K\"ahler-Einstein metric on $X$, but, as this happens,
a substantial amount of 
trace-free Ricci curvature accumulates in  the  neck joining the two copies of $X$. 
 \end{xpl}

The same techniques used to prove this result also have a number
of other geometrical ramifications. Of these, we will  highlight
the following:

\begin{main} 
Let $X_j$, $j=1, \ldots, 4$ be  smooth, compact 
almost-complex $4$-manifolds for which the 
 mod-2 Seiberg-Witten invariant is non-zero, and 
suppose, moreover, that  
(\ref{tom}--\ref{harry}) hold.  Let $N$ be any oriented $4$-manifold
with $b_+=0$.  Then, for any $m=2, 3$ or    $4$,  the smooth $4$-manifold 
$[ \#_{j=1}^mX_j]\#N$ does not admit 
Einstein metrics if 
$$
 12(m-1) +   (12 b_1 + 3b_- )(N)    \geq \sum_{j=1}^m c_1^2(X_j) .
$$
\end{main}

This generalizes the key technical result of  \cite{il1}.

\section{The Seiberg-Witten Equations}\label{curv}

If $M$ is a smooth  oriented $4$-manifold, there always exist 
Hermitian line bundles $L\to M$ such that $c_{1}(L)\equiv w_{2}(TM)
\bmod 2$. For any such $L$, and any Riemannian metric $g$ on $M$,
one can then find rank-$2$ Hermitian vector bundles ${\mathbb V}_\pm$
which formally satisfy
$${\mathbb V}_\pm= {\mathbb S}_{\pm}\otimes L^{1/2},$$
where ${\mathbb S}_{\pm}$ are the locally defined left- and
right-handed spinor bundles of $(M,g)$. 
Such a choice of   ${\mathbb V}_{\pm}$, up to isomorphism, 
is called a 
spin$^{c}$ structure $\mathfrak c$ on $M$, and is determined, modulo
the $2$-torsion subgroup of $H_{1}(M,\ZZ)$, by the first Chern class
$c_{1}(L)= c_{1}({\mathbb V}_{\pm})
\in H^{2}(M,\ZZ )$ of the spin$^{c}$ structure. 
Every unitary 
connection $A$ on $L$ induces a Dirac operator
$$D_{A}: \Gamma ({\mathbb V}_{+})\to \Gamma ({\mathbb V}_{-}).$$
If $A$ is such a connection, and if $\Phi$ is a section of
${\mathbb V}_{+}$, the pair $(\Phi , A)$ is said to satisfy the
 {\em Seiberg-Witten equations} \cite{witten} if  
\begin{eqnarray} D_{A}\Phi &=&0\label{drc}\\
 F_{A}^+&=&i\sigma (\Phi ),\label{sd}\end{eqnarray}
were 
$F_{A}^+$  is the self-dual part of
the curvature of $A$, and where $\sigma : {\mathbb V}_+\to \Lambda^+$
is a real-quadratic map arising from the isomorphism 
$\Lambda^{+}\otimes \CC = \odot^{2}{\mathbb S}_{+}$.
 
For the $4$-manifolds of primary interest here, there turn out to be 
certain spin$^{c}$ structures  for which 
there exists a solution of the Seiberg-Witten equations 
for each metric $g$. This situation is neatly codified by the 
following 
 terminology, first introduced by Kronheimer \cite{K}:
 
\begin{defn}
Let $M$ be a smooth compact oriented $4$-manifold
with $b_{+}\geq 2$. An element $a\in  H^{2}(M,\ZZ )/
\mbox{\rm torsion}$ will be called a {\bf monopole
class} of $M$ iff there exists a spin$^{c}$ structure
$\mathfrak c$ 
on $M$ with 
$$c_{1}^{\RR}(L)= a$$ which has the property
that  the corresponding  Seiberg-Witten 
equations (\ref{drc}--\ref{sd})
have a solution for every Riemannian  metric $g$ on $M$. 
Here $c_{1}^{\RR}(L)$ denotes the image of the first Chern class $c_1(L)$ 
of the anti-canonical line bundle of $\mathfrak c$ in 
$H^{2}(M,\ZZ )/
\mbox{\rm torsion}\subset H^2(M, \RR)$. 
\end{defn}

As was first pointed out by Witten \cite{witten},
the existence of a monopole class implies an  
{\em a priori} lower bound on the $L^{2}$ norm of 
the scalar curvature of Riemannian metrics. 
The sharp form of this  estimate reads as follows \cite{lpm}: 

 \begin{prop} \label{west}
Let $M$ be a   smooth compact  
oriented  $4$-manifold with monopole class 
$a\in H^{2}(M,\ZZ)/\mbox{\rm torsion}\subset H^{2}(M,\RR )$. 
 If $a\neq 0$, 
then $M$ does not admit metrics of positive scalar curvature. 
Moreover, if 
 $g$ is any Riemannian metric on $M$, and 
if $a^{+}\in H^{2}(M,\RR )$ denotes the 
self-dual part of $a$ with respect to the 
 decomposition 
$${\mathcal H}^{2}_{g}= {\mathcal H}^{+}_{g}\oplus {\mathcal H}^{-}_{g}$$
of  the
space of $g$-harmonic $2$-forms into eigenspaces of the 
$\star$ operator, 
then the scalar curvature $s$  of $g$ satisfies 
$$
\int_{M}  s^{2} d\mu \geq 32\pi^2  
(a^{+})^{2} ,
$$
where $d\mu$ denotes the Riemannian volume form of 
$g$.   Finally,  the  inequality 
is always strict unless  $M$ admits a  complex structure of K\"ahler type,  
with first Chern class $a$. 
\end{prop}

One important formal consequence is the following; cf. \cite{il1}.

\begin{prop}
 Let $M$ be any smooth compact oriented $4$-manifold with 
 $b_+(M)\geq 2$. Then  ${\mathcal C}=\{ \mbox{monopole classes of }M\}$ 
  is a finite set. 
 \end{prop}
 \begin{proof}
 Let $g_1=g$ be any Riemannian metric on $M$, and let 
 ${\mathbf e}_1=[\alpha_1]$ be the cohomology class of a harmonic self-dual
 form with respect to $g$, normalized so that ${\mathbf e}_1^2:=[\alpha_1]\cdot [\alpha_1]=1$.
Now the local surjectivity of the period map \cite[Prop. 4.3.14]{don}
asserts that $[\alpha_1]$ has a neighbohood in $H^2(M,\RR )$
in which every element can be represented by a self-dual harmonic form
for some perturbation of $g$. But any open set in a finite-dimensional
vector space spans the space. Thus we can find a basis $\{ {\mathbf e}_j= [\alpha_j]~|~j= 1, \ldots ,
b_2 (M)\}$ for $H^2(M,\RR )$ and a collection of Riemannian metrics $\{ g_j ~|~j= 1, \ldots ,
b_2 (M)\}$ such that the harmonic representative $\alpha_j$ of ${\mathbf e}_j$ 
with respect to $g_j$ is self-dual.  For convenience, we  normalize these basis
elements so that ${\mathbf e}_j^2=1$ for each $j$. Let $L_j: H^2(M,\RR ) \to \RR$
denote the linear functionals $L_j(x)= {\mathbf e}_j\cdot x$, and now notice that
Proposition 
\ref{west}, together with 
the Cauchy-Schwarz inequality, tells  us that any monopole class 
$a\in H^2(M, \ZZ)/\mbox{torsion}$ satisfies 
$$|L_j (a)| = |[\alpha_j ] \cdot a^+_{g_j}|\leq \sqrt{(a^+_{g_j})^2} \leq 
\left(\int_M s_{g_j}^2d\mu_{g_j}\right)^{1/2}=\varrho_j,$$
where the constant $\varrho_j$ is independent of $a$. This shows that 
 ${\mathcal C}\subset H^2(M,\RR )$
is contained in the $b_2(M)$-dimensional parallelepiped 
$$\left\{ x\in H^2 (M,\RR ) ~\Big|~ |L_j(x)|\leq \varrho_j ~\forall j= 1,\ldots , b_2(M)\right\},$$
which, while geometrically wide and flat, is nonetheless a compact set. 
Since ${\mathcal C}\subset H^2(M,\ZZ )/\mbox{torsion}$ is
also discrete, it follows  that $\mathcal C$ is finite. 
 \end{proof}

More
recently, the second author  discovered  that equations 
(\ref{drc}--\ref{sd}) also imply a family of
 estimates involving the 
self-dual Weyl curvature  \cite{lric}:

 \begin{prop}[LeBrun]\label{est}
Let $M$ be a   smooth compact  
oriented  $4$-manifold with monopole class 
$a$, and let 
 $g$ be any Riemannian metric on $M$. 
Then the scalar curvature $s$ and self-dual Weyl curvature
$W_{+}$ of $g$ satisfy 
$$
\int_{M} (s- \sqrt{6}  
|W_{+}|)^{2} d\mu \geq 72\pi^2  
(a^{+})^{2} ,
$$
where $d\mu$ denotes the Riemannian volume form of 
$g$, and where the point-wise norms are calculated with
respect to $g$. Moreover,  the inequality is always strict
unless  $M$ admits a 
 symplectic structure with  first Chern class $a$.
 \end{prop}

We remark in passing that one can actually say a great deal more about
the borderline cases of the above inequalities. 
Indeed, the inequality of  Proposition \ref{west}
is only saturated by constant-scalar-curvature K\"ahler metrics \cite{lpm}, whereas the 
inequality of Proposition \ref{est} is only saturated by a class of 
{almost}-K\"ahler metrics with  extremely special curvature properties
\cite{lric}. 
However, the weaker statements given above will
be quite sufficient for the purposes of the present article.

\section{The Bauer-Furuta Invariant}\label{bauer}

Bauer and Furuta \cite{baufu,bauer2} recently introduced a remarkable new 
method for detecting monopole classes. The crux of their theory is the   definition of  a 
generalized Seiberg-Witten 
invariant which takes values in an   
 equivariant stable cohomotopy group. As a courtesy to the reader, we therefore 
begin with a brief introduction to
the relevant homotopy-theoretic machinery. 

If $({\mathfrak X},x)$ is any finite cell complex with base point, its {\em suspension} $\Sigma {\mathfrak X}=\Sigma^1 {\mathfrak X}$ 
is just the one-point compactification of 
$({\mathfrak X}-\{ x\}) \times \RR$, equipped with  base point $\infty$. Setting
$\Sigma^{n}{\mathfrak X}= \Sigma (\Sigma^{n-1} {\mathfrak X})$ for  $n> 1$, one immediately sees that 
 $\Sigma^n {\mathfrak X}$ is similarly the one-point compactification of 
$({\mathfrak X}-\{ x\}) \times \RR^n$. 
 Now if  $({\mathfrak Y},y)$ is any
other pointed finite cell complex, then the set $[\Sigma^n{\mathfrak X} , {\mathfrak Y}]$ of based homotopy classes of 
 maps $\Sigma^n{\mathfrak X}\to {\mathfrak Y}$ can be identified with $\pi_n (\{ \mbox{based maps } {\mathfrak X}\to {\mathfrak Y}\})$,
and so is an Abelian  group if $n\geq 2$. There is a natural 
homomorphism  $[\Sigma^n {\mathfrak X} , \Sigma^n {\mathfrak Y}]\to [\Sigma^{n+1}{\mathfrak X}, \Sigma^{n+1}{\mathfrak Y}]$,
induced by Cartesian products with the identity $\RR \to \RR$. If $n$ is large with 
respect to the dimensions of ${\mathfrak X}$ and ${\mathfrak Y}$, then this homomorphism turns out to be
an isomorphism. In particular, we may define the $m^{th}$ {\em stable cohomotopy
group} of ${\mathfrak X}$ by
$$\pi^m ({\mathfrak X}) = \lim_{n\to \infty} [\Sigma^n {\mathfrak X} , \Sigma^{n}S^m ] = \lim_{n\to \infty}
[\Sigma^n {\mathfrak X} , S^{m+n} ].$$

The Bauer-Furuta construction  also involves  $S^1$-equivariant 
analogs of these groups. Suppose that ${\mathfrak X}$ is equipped with an $S^1$-action
leaving the base point $x$ fixed. Then we may equip $\Sigma^{3n}{\mathfrak X}$
with an $S^1$-action by choosing to view it as the one-point
compactification of $({\mathfrak X}-\{ x\}) \times \CC^n \times \RR^n$, and letting 
$S^1\subset \CC^\times$ act on the $\CC^n$ factor by scalar multiplication. If $({\mathfrak Y},y)$
is another complex with $S^1$-action, then the set $[\Sigma^{3n}{\mathfrak X}, {\mathfrak Y}]_{S^1}$
of based homotopy classes of $S^1$-equivariant maps $\Sigma^{3n}{\mathfrak X} \to {\mathfrak Y}$
can be identified with $\pi_n (  \{\mbox{based $S^1$-equivariant maps  } \Sigma^{2n}{\mathfrak X}\to {\mathfrak Y}\})$,
and so is an Abelian group whenever $n \geq 2$. We again have 
natural homomorphisms $[\Sigma^{3n} {\mathfrak X} , \Sigma^{3n} {\mathfrak Y}]_{S^1}\to [\Sigma^{3(n+1)}{\mathfrak X}, \Sigma^{3(n+1)}{\mathfrak Y}]_{S^1}$, and these homomorphisms are again isomorphisms 
for $n$ sufficiently large.  The $S^1$-equivariant cohomotopy groups of ${\mathfrak X}$ 
used in the  Bauer-Furuta construction 
are then given  by 
$$\pi^m_{S^1} ({\mathfrak X}) = \lim_{n\to \infty} [\Sigma^{3n} {\mathfrak X} , \Sigma^{3n}S^m ]_{S^1} = \lim_{n\to \infty}
[\Sigma^{3n} {\mathfrak X} , S^{m+3n} ]_{S^1},$$
where $S^1$ acts on ${\mathfrak X}$ in the given manner, and acts trivially on $S^m$.
In particular, if ${\mathbb E}^+= {\mathbb E}\cup \{\infty \}$ is the Thom space of a 
complex vector bundle ${\mathbb E}\to B$ over a compact base, $S^1$ will be understood to act on 
${\mathbb E}$ by fiber-wise complex multiplication, and
$$\pi^m_{S^1} ({\mathbb E}^+)=\pi_0 ( \left\{\mbox{$S^1$-equivariant  maps } 
({\mathbb E}\oplus \CC^n \oplus \RR^{n})^+ \to  (\CC^n \times \RR^{m+n})\cup \{\infty\}\right\})$$
for $n$ sufficiently large. 
Notice that this only depends on the 
stable equivalence class of ${\mathbb E}$, rather than ${\mathbb E}$ itself. In particular, 
if ${\mathbb E}$ is merely specified by an element of $K^0 (B)$, we can still
define the group in question by representing the virtual bundle ${\mathbb E}$ as the  
formal difference ${\mathbb F}\ominus \CC^\ell$ of a complex vector bundle 
and a trivial bundle, and then taking 
$$\pi^m_{S^1} ({\mathbb E}^+)=\pi_0 ( \left\{\mbox{$S^1$-equivariant  maps } 
({\mathbb F}\oplus \CC^{n-\ell}\oplus \RR^{n})^+ \to (\CC^n \times \RR^{m+n})
\cup \{\infty\}\right\})$$
for $n$ sufficiently large. 

Bauer and Furuta \cite{baufu}  observe  that one can associate stable cohomotopy 
classes with  sufficiently nice maps between Hilbert spaces. Indeed, let
$${\mathcal L} :  {\mathcal H}_1 \to {\mathcal H}_2$$
be a Fredholm linear map between real Hilbert spaces, and 
let 
$${\mathcal F} :  {\mathcal H}_1 \to {\mathcal H}_2$$
be a continuous map of the form 
$${\mathcal F} = {\mathcal L} + {\mathcal K},$$
where $\mathcal K$ is a compact  (non-linear) operator.
Assume, in addition, that ${\mathcal F}$ satisfies the following 
{\bf boundedness condition}: 
\begin{center}
{\sl the inverse image of any bounded set is bounded.} 
\end{center}
The boundedness condition allows us to extend 
${\mathcal F}$ to a continuous map  
$${\mathcal F}^+={\mathcal H}_1^+ \to {\mathcal H}_2^+$$
where ${\mathcal H}_j^+= {\mathcal H}_j\cup \{ \infty\}$
is identified with the unit sphere in $S({\mathcal H}_j\oplus \RR)
\subset {\mathcal H}_j\oplus \RR$
via stereographic projection.

Now given any  finite-dimensional subspace  ${\mathbf W}\subset {\mathcal H}_2$ 
with  $\mbox{coker }{\mathcal L} \subset {\mathbf W}$, consider corresponding  
finite-dimensional subspace ${\mathbf V}:= {\mathcal L} ^{-1}({\mathbf W}) 
\subset {\mathcal H}_1$. We will then  say that ${\mathbf W}$ is
a {\em good} subspace if ${\mathcal F}({\mathbf V})$ does not intersect the 
unit sphere $S({{\mathbf W}^\perp})\subset {\mathbf W}^\perp$.  For any good subspace ${\mathbf W}$ we thus 
have  a map $\mbox{pr}_{\mathbf W}\circ {\mathcal F}^+|_{{\mathbf V}^+} : {\mathbf V}^+\to {\mathbf W}^+$ between  
finite-dimensional spheres, where the deformation retraction 
$\mbox{pr}_{\mathbf W}: {\mathcal H}_2^+ - 
S({{\mathbf W}^\perp}) \to {\mathbf W}^+$ corresponds, via stereographic projection,  to the 
map $S({\mathcal H}_2 \oplus \RR)-S({\mathbf W}^\perp) \to S({\mathbf W}\oplus \RR)$
given by $v\to v^\|/|v^\| |$, where  $v^\|$ is the orthogonal projection of
$v$ into ${\mathbf W}\oplus \RR$. It
 turns out that any finite-dimensional subspace 
of ${\mathcal H}_2$ is a proper subspace of some good subspace,
and we can thus define a stable cohomotopy class $[{\mathcal F}]
\in \pi^{-\mbox{ind} ({\mathcal L})}(S^0)$ by taking the limit of  this
construction over  a nested sequence  ${\mathbf W}_1\subset {\mathbf W}_2\subset {\mathbf W}_3
\subset \cdots$ of good subspaces.
In fact, the   cohomotopy classes constructed from the individual  ${\mathbf W}_j$ 
all agree with one another  for sufficiently large $j$. Moreover,
$[{\mathcal F}]$  is independent of   all  the choices
made in its construction.

The equivariant version of the construction is similar. Suppose that 
 ${\mathcal L}$ is the direct sum ${\mathfrak L}_\CC\oplus {\mathfrak L}_\RR$ of 
a complex-linear Fredholm map ${\mathfrak L}_\CC$ between complex Hilbert spaces and
a real linear Fredholm map ${\mathfrak L}_\RR$ between real Hilbert spaces. 
Also assume that $\mathcal F$ is an $S^1$-equivariant map. Then 
the same procedure gives us an equivariant cohomotopy class
$[{\mathcal F}]
\in \pi^{\jmath}_{S^1}(\CC^\imath\cup \{\infty\})$, where 
 $\imath= \mbox{ind}_\CC {\mathfrak L}_\CC$ and $\jmath= -\mbox{ind}_\RR {\mathfrak L}_\RR$.

If we replace ${\mathcal H}_1$ and ${\mathcal H}_2$ with real Hilbert-space bundles over a
finite cell complex $B$,  and allow $\mathcal L$ to be a 
Fredholm morphism of Hilbert bundles, then much the same construction 
still applies. Indeed, we can always find a trivialization
${\mathcal H}_2\cong B \times {\mathcal H}$ so that there is a fixed finite-dimensional
subspace  ${\mathbf W}_0\subset {\mathcal H}$ with $\mbox{coker } {\mathcal L}_x
\subset {\mathbf W}_0$ for every $x\in B$. A finite-dimensional subspace
${\mathbf W}\supset {\mathbf W}_0$ is then said to be good if 
${\mathcal F}( {\mathbf V} )$ avoids the unit sphere in  ${\mathbf W}^\perp$,
where ${\mathbf V} ={\mathcal L}^{-1}({\mathbf W})$ is now a bundle. 
For each good subspace $\bf W$,  the map 
$\mbox{pr}_{\mathbf W}\circ {\mathcal F}^+|_{{\mathbf V}^+} : {\mathbf V}^+\to {\mathbf W}^+$
determines an element of   $\pi^0 ( {\mathbb E}^+)$, where 
${\mathbb E}={\mathbf V} \ominus {\mathbf W}$ is the virtual index of ${\mathcal L}$ in $KO^0 (B)$,
and ${\mathbb E}^+$ is its
Thom space.  It once again turns out that 
 any finite-dimensional subspace of $\mathcal H$ is a proper
subspace of a good subspace, so we may take the limit over a nested sequence 
 ${\mathbf W}_1\subset {\mathbf W}_2\subset {\mathbf W}_3
\subset \cdots$ of good subspaces to define the  stable cohomotopy  class
 $[\mathcal F]\in \pi^0 ( {\mathbb E}^+)$, and one can 
 then show that this class  is independent of 
all choices made along the way. 

Similarly, there is an $S^1$-equivariant version of this 
 story. The relevant special case goes as follows: suppose
that ${\mathcal H}_1$ and ${\mathcal H}_2$ are each direct sums of 
complex Hilbert-space bundles over $B$ and {\em trivial} real Hilbert-space bundles.
Suppose also that  ${\mathcal L}$ decomposes as ${\mathfrak L}_\CC\oplus {\mathfrak L}_\RR$,
where ${\mathfrak L}_\CC$ is a complex-linear Fredholm morphism, and 
${\mathfrak L}_\RR$ is a Fredholm linear map between real Hilbert spaces.
Then if ${\mathcal F} = {\mathcal L}+ {\mathcal K}$ is $S^1$-equivariant and 
satisfies the boundedness condition, with ${\mathcal K}$ is compact, 
then the previous construction, carried out $S^1$-equivariantly, gives us a
stable cohomotopy class $[{\mathcal F}]\in \pi^\jmath_{S^1} ({\mathbb E}^+)$, 
where ${\mathbb E}\in K^0(B)$ is the virtual index bundle of $\mathcal L$,
${\mathbb E}^+$ is its Thom space, and  $\jmath= -\mbox{ind}_\RR {\mathfrak L}_\RR$.

We now describe how Bauer and Furuta use this machinery to construct
$4$-manifold invariants. 
Let $(M,g)$ be a smooth compact oriented Riemannian $4$-manifold with 
a chosen spin$^c$ structure $\mathfrak c$. Choose some 
smooth $U(1)$ connection $A_0$ on 
the anti-canonical line bundle $L$ of $\mathfrak c$,
and, for some $k> 2$, consider the {\em monopole  map} 
$$\begin{array}{rclcrccccl}
L^2_{k+1} ({\mathbb V}_+) & \oplus & L^2_{k+1} (\Lambda^1)
&\stackrel{{\psi}_k}{\longrightarrow}& 
L^2_k({\mathbb V}_-)& \oplus& L^2_k(\Lambda^+) &\oplus&L^2_k/\RR
\\ &&&&&&&\\
\Big( ~~ \Phi & , &\alpha ~~\Big) & \mapsto & \Big( ~~
D_{A}\Phi &, &iF_A^+ + \sigma (\Phi ) & 
,&
d^\ast \alpha ~~ \Big)
\end{array}$$
of real Hilbert spaces. Here $A:=A_0+i\alpha$,  and $\sigma : {\mathbb V}_+\to \Lambda^+$
is the quadratic map occurring in  equation (\ref{sd}).  
The monopole map really just represents  the Seiberg-Witten equations (\ref{drc}--\ref{sd})
and a  gauge-fixing condition, insofar as 
 $\psi_k^{-1}(0)$ by definition consists of the $L^2_{k+1}$ solutions of the 
Seiberg-Witten equations  for which  
$d^*(A-A_0)=0$. 
Now notice that 
$$\psi_k = {\mathcal L} + {\mathcal K},$$
where ${\mathcal L}$ is the linear Fredholm  map 
given by 
$${\mathcal L}(\Phi, \alpha) = (D_{A_0}\Phi , -d^+\alpha , d^*\alpha )$$
and where ${\mathcal K}$ is the compact non-linear operator
$${\mathcal K}(\Phi, \alpha) = (\frac{i}{2} \alpha \cdot \Phi , iF^+_{A_0}+\sigma (\Phi) , 0 ).$$
Also notice that  $\psi_k$ is equivariant with respect to the $S^1$ action 
associated with the complex-Hilbert-space structures of $L^2_{k+1}({\mathbb V}_+)$
and $L^2_{k}({\mathbb V}_-)$. 
If $b_1(M)=0$, and $k>3$, a Weitzenb\"ock-formula argument \cite{baufu}
shows that $\psi_k$ satisfies the boundedness condition. By the previous
discussion,   it therefore determines an equivariant
stable cohomotopy
class $[\psi_k ]\in \pi^{b}_{S^1}(\CC^\imath\cup \{\infty\})$, where $b=b_+(M)$ and 
$\imath=[c_1^2(L)-\tau (M)]/8$. One can  show that 
 this class  is independent of the metric $g$ used to 
construct it, and so constitutes a smooth manifold  invariant of $M$, depending
only on the spin$^c$ structure $\mathfrak c$. It  is also independent of 
the chosen integer $k > 3$, so we will therefore usually denote the
invariant by $[\psi]$, rather than by $[\psi_k]$. For much the same reason, 
we will simply call the 
monopole map   $\psi$, rather than $\psi_k$, 
whenever this is unlikely to lead to confusion. 

When $b_1(M)\neq 0$, the story is somewhat more involved, but the 
upshot is similar. Choose a base-point $\ast \in M$, and  identify $H^1(M,\ZZ )$ with the group of 
{\em harmonic maps} $u: M\to S^1$ such that  $u(\ast )=1$.
Then  ${\psi}_k$ is   equivariant  with respect to the $H^1(M,\ZZ )$-action by 
gauge transformations
\begin{eqnarray*}
(\Phi , \alpha) & \mapsto & (u\Phi, \alpha -2i ~d\log u)\\
(\Psi , \varphi , c, f) & \mapsto & (u\Psi, \varphi , c , f)
\end{eqnarray*}
and we may use these actions to   define 
\begin{eqnarray*}
{\mathcal A}_{k+1}&=& \Big[
L^2_{k+1} ({\mathbb V}_+)  \oplus L^2_{k+1} (\Lambda^1)\Big]/H^1(M,\ZZ )\\
{\mathcal B}_k&=& 
  \Big[\Big(L^2_k({\mathbb V}_-) \oplus L^2_k(\Lambda^+) \oplus L^2_k/\RR\Big)
\times H^1(M,\RR )\Big]/H^1(M,\ZZ ).
\end{eqnarray*}
Observe that both ${\mathcal A}_{k+1}$ and ${\mathcal B}_k$ are
 Hilbert-space bundles over 
the Picard torus 
$$\mbox{Pic}^0(M)= H^1(M, \RR)/H^1(M, \ZZ ), $$
where the projection  ${\mathcal A}_{k+1}\to Pic^0(M)$ sends
$(\Phi , \alpha )$ to the  
 harmonic part of $\alpha$. 
 The monopole map now becomes an $S^1$-equivariant  compact perturbation 
 of a Fredholm morphism of Hilbert bundles, and we may associated to it
 a stable cohomotopy class in the manner described above. 
 That is, by 
choosing
a suitable trivialization of ${\mathcal B}_k$, we may view our monopole map as an
$S^1$-equivariant 
smooth map 
$$\psi_k : {\mathcal A}_{k+1}\to {\mathcal C}_{k},$$
where
$${\mathcal C}_k= L^2_k({\mathbb V}_-) \oplus L^2_k(\Lambda^+) \oplus L^2_k/\RR ,$$
and extract its stable cohomotopy class by restricting $\psi_k$ to 
${\mathbf V}= {\mathcal L}^{-1}({\mathbf W})$ and then projecting into 
${\mathbf W}^+$ for a  good subspace ${\mathbf W}\subset {\mathcal C}_k$
of sufficiently large dimension.
This  equivariant stable cohomotopy class  is  independent of $g$ and $k$, 
and the resulting picture may be codified as follows:

\begin{thm}[Bauer-Furuta]
Let $M$ be a smooth compact oriented $4$-manifold with 
$b_+(M)=b$, and let $\mathfrak c$ be a 
spin$^c$ structure on $M$. 
Then the corresponding  monopole map $\psi$ determines a stable cohomotopy class 
$$[\psi] \in \pi^b_{S^1} (\mbox{ind} (D)^+), 
$$
where $ind (D)\in K^0(\mbox{Pic}^0 (M))$ is the virtual 
index bundle of the spin$^c$ Dirac operator.
Moreover, if $b_+(X) \geq b_1(X) + 2$, an orientation of the vector space 
$H^1(M,\RR )\oplus {\mathcal H}_g^+$ determines   a 
group homomorphism $ \pi^b_{S^1} (\mbox{ind} (D)^+)
\to {\mathbb Z}$ which maps $[\psi ]$ to the usual Seiberg-Witten invariant. 
\end{thm}

The class $[\psi]$ will henceforth be called 
the  {\em  Bauer-Furuta invariant} of $(M, {\mathfrak c})$. 
Its importance for our purposes stems from the following observation:

\begin{prop} \label{woosh}
If the Bauer-Furuta invariant $[\psi]$ of $(M, {\mathfrak c})$ is 
non-zero, then  $c_1^{\RR} (L)$ is a monopole class. (Here, once again, 
$L$ is the anti-canonical line bundle of $\mathfrak c$, and $c_1^{\RR} (L)$
is the image of $c_1(L)$ in $H^2(M, \ZZ)/\mbox{torsion}$.) 
\end{prop}
\begin{proof}
Given a smooth metric $g$ on $M$, we will  show that the  Seiberg-Witten 
equations (\ref{drc}--\ref{sd}) admit a solution $(\Phi, A)$ for the spin$^c$
structure $\mathfrak c$. Our argument will hinge on the fact that, while  the monopole
map $\psi_k$ used in the construction  depends on the choice of an  integer $k > 3$,  we always have a
commutative diagram
\begin{eqnarray*}
{\mathcal A}_{k} & \stackrel{\psi_{k-1}}{\longrightarrow} & {\mathcal C}_{k-1}
 \\
\uparrow &  & \uparrow \\
{\mathcal A}_{k+1} & \stackrel{\psi_k}{\longrightarrow} & {\mathcal C}_k
\end{eqnarray*}
where $\psi_{k-1}$ is also  continuous. 
The result will be deduced  from this by
invoking the Rellich theorem. 

Indeed, consider the bounded 
self-adjoint operator $Q: {\mathcal C}_k \to {\mathcal C}_k$
defined by 
$$\langle v , Q u\rangle_{L^2_k} := \langle v , u \rangle_{L^2_{k-1}}.$$
Because the positive self-adjoint operator $Q$ is   compact by the Rellich theorem, its
 eigenspaces ${\mathbf E}_j$ are  mutually orthogonal, finite-dimensional, 
span
 ${\mathcal C}_k$, and may be taken to be  listed in an order such  that  
 the   corresponding sequence of  eigenvalues $\lambda_1 > \lambda_2 > \cdots > \lambda_j
> \cdots  > 0$
 monotonically decreases to $0$.  Now choose a sequence of
good subspaces ${\mathbf W}_1\subset {\mathbf W}_2 \subset {\mathbf W}_3 \subset \cdots $ for 
$\psi_k$ so that $\oplus_{i=1}^j{\mathbf E}_i\subset {\mathbf W}_j$. 
If, for large $j$, $\psi_k ({\mathbf V}_j)$ were to avoid the unit ball in ${\mathbf W}_j^\perp$,
then $(\mbox{pr}_{{\mathbf W}_j} \circ \psi_k) ({\mathbf V}_j)$ would avoid the origin, and
 $\mbox{pr}_{{\mathbf W}_j} \circ \psi_k: {\mathbf V}_j\cup \{ \infty \}\to {\mathbf W}_j \cup  \{ \infty \}$
would be $S^1$-equivariantly homotopic to the constant map at $\infty$
via dilations; but we would then have $[\psi]=0$, contradicting our hypothesis. 
Hence the unit ball in ${\mathbf W}_j^\perp$ must meet the image of $\psi_k$ for
 all j. That is, there must be  sequences $p_j\in {\mathcal A}_{k+1}$,
$q_j\in  {\mathcal C}_k$, with 
$$\psi_k (p_j) = q_j$$
and 
$$\|q_j\|_{L^2_k} < 1, ~~ q_j \perp \oplus_{i=1}^j {\mathbf E}_i.$$
However, 
$$u \perp \oplus_{i=1}^j{\mathbf E}_i ~~\Longrightarrow ~~\|u\|^2_{L^2_{k-1}}
\leq \lambda_{j+1}\|u\|^2_{L^2_k},$$
so it  follows that 
that $\| q_j\|_{L^2_{k-1}}\leq \sqrt{\lambda_{j+1}}\to  0,$ and hence $q_j \to 0$ in $L^2_{k-1}$.

On the other hand, because  $\psi_k$ satisfies the boundedness condition, and because 
the $q_j$ are all in the unit ball, we have 
$$\|p_j\|_{L^2_{k+1}} < C$$
for some $C$. The Rellich theorem therefore tells us that there is a subsequence 
$p_{j_i}$ which converges in the $L^2_k$ norm to some $p_\infty \in {\mathcal A}_{k}$.
 The continuity of $\psi_{k-1}$, together with its compatibility with $\psi_k$, now tells us that 
$$
\psi_{k-1} (p_{\infty})  =  \psi_{k-1} (\lim_{i\to \infty} p_{j_i})
 =\lim_{i\to \infty} \psi_{k-1} ( p_{j_i})=  \lim_{i\to \infty} q_{j_i} =0 ,
$$
which is to say 
 that $p_\infty= (\Phi, \alpha )$ 
gives us 
 a solution $ (\Phi, A_0 + i\alpha )$ of the gauge-fixed Seiberg-Witten equations 
of class $L^2_k$. The usual bootstrap regularity argument then shows that
this solution is actually  $C^\infty$. 
 \end{proof}

One of Bauer's key results involves    
the {\em smash product} $$\wedge : \pi^m_{S^1}({\mathfrak X}) \times 
 \pi^{m'}_{S^1}({\mathfrak X}')\to \pi^{m+m'}_{S^1}({\mathfrak X}\wedge {\mathfrak X}').$$   Here, if $({\mathfrak X},x)$ and $({\mathfrak X}',x')$ are 
pointed finite complexes, their {\em join} 
${\mathfrak X}\wedge {\mathfrak X}'$ is just the one-point compactification of 
$({\mathfrak X}-x)\times ({\mathfrak X}'-x')$. The operation in question is then induced by 
taking Cartesian products. Namely, if 
$$f : ({\mathfrak X}-x)\times \CC^n \times \RR^n \to \CC^n \times \RR^{m+n}$$
and
$$f' : ({\mathfrak X}'-x')\times \CC^{n'} \times \RR^{n'} \to\CC^{n'} \times \RR^{m'+n'}$$
are $S^1$-equivariant  proper maps, so is
$$f \times f' :   ({\mathfrak X}-x)\times ({\mathfrak X}'-x') \times \CC^{n+n'} \times \RR^{n+n'}
\to  \CC^{n+n'} \times \RR^{m+m'+n+n'},$$
and setting 
$[f] \wedge [f']= [f \times f']$ then gives us the desired operation.

\begin{thm}[Bauer] \label{smash}
Let $X$ and $Y$ be smooth compact  oriented 4-manifolds
with  spin$^c$ structures ${\mathfrak c}_X$ and ${\mathfrak c}_Y$.
Then the Bauer-Furuta invariant of the connected sum 
 $X\# Y$, equipped with the spin$^c$ structure ${\mathfrak c}_X\# {\mathfrak c}_Y$,
  is exactly the  the smash product of the 
invariants of $(X, {\mathfrak c}_X)$ and $(Y,{\mathfrak c}_Y)$:
$$[\psi_{X\# Y}] = [\psi_{X}]  \wedge [\psi_{Y}] .$$
In other words, $\psi_{X\# Y}$ and $\psi_{X}\times \psi_{Y}$
represent the same equivariant 
stable cohomotopy class. 
\end{thm}

\begin{cor}\label{smoosh}
Let $M$ be a  smooth compact  oriented 4-manifold with 
a ${spin}^{c}$-structure
${\mathfrak c}_{M}$  for which the 
Bauer-Furuta invariant is 
non-trivial. Let $N$ be a smooth compact  oriented 4-manifold  with
negative definite
 intersection form, and let ${\mathfrak c}_{N}$ be any  
 ${spin}^{c}$-structure on $N$ with $c_1^2=-b_2(N)$. Then 
the Bauer-Furuta invariant is also non-trivial for the 
spin$^c$ structure ${\mathfrak c}_{M} \# {\mathfrak c}_{N}$
on the connected sum $M \# {N}$. 
\end{cor}

\begin{proof} Since we have chosen a spin$^c$ structure on
$N$ for which the Dirac operator has index $[c_1^2-\tau(N)]/8 =0$, 
the vector  bundle $\mbox{ind} (D_N)\to \mbox{Pic}^0(N)$
has rank $0$, and the image of the Bauer-Furuta invariant of $N$
under the map 
$$\pi^0_{S^1}(\mbox{ind} (D_N)^+) \to  \pi^0_{S^1}(\mbox{pt}^+),$$
induced by restriction to the fiber over a point in $\mbox{Pic}^0(N)$,  is   just 
 the class of the identity map.
By Theorem \ref{smash}, the restriction of the invariant 
$[\psi ] = [\psi_M \wedge \psi_N]$  of $M\#N$ to 
$\pi^b_{S^1}( [ \mbox{ind}(D_{M\# N})|_{\mbox{Pic}^0(M) \times \mbox{pt}}]^+)$ 
 therefore equals  $[\psi_M] \wedge[ \mbox{id} ]= [\psi_M]$.
\end{proof} 

When $b_1(M)=0$ and the spin$^c$ structure arises from an almost-complex structure,
one has 
$$ 
\pi^b_{S^1} (\mbox{ind} (D)^+ )=\pi^b_{S^1} ( S^{b+1}),$$
where the sphere $S^{b+1}$ is thought of as 
$\CC^\imath
\cup\{ \infty\}$
for  $\imath=[c_{1}^{2}(L)-\tau(X)]/8= (b+1)/2$.   
If the mod-$2$ Seiberg-Witten invariant is non-zero and $b=b_+(M)\equiv 3 \bmod 4$,
Bauer  \cite{bauer2} then  shows that the Bauer-Furuta invariant is  represented by 
the suspension of the Hopf map $S^3 \to S^2$. 
In conjunction with  Theorem  \ref{smash}, 
this then implies the following:

\begin{thm}[Bauer] \label{dimsum}
For each $j\in \{ 1, 2, 3, 4\}$,
	let $(X_j,J_j)$    be a smooth  compact  almost-complex 
	$4$-manifold with $b_{1}=0$ and  $b_{+}\equiv 3\bmod 4$.
         Suppose, moreover, that $\sum_{j=1}^4b_{+}(X_{j})\equiv 4 \bmod 8$,
         and that  
	  the canonical spin$^{c}$ structure on each $(X_j.J_j)$ 
          has non-zero 
         mod-$2$ Seiberg-Witten.  Then the corresponding 
	spin$^{c}$ structures on the connected sums  $\#_{j=1}^mX_j$, $m=1, \ldots, 4$, 
          all   have non-zero Bauer-Furuta 	invariant.
	\end{thm}

Indeed, when $m =2$ or $3$, this can  be deduced  by just looking at the 
image of the Bauer-Furuta invariant in the non-equivariant  stable homotopy group 
$\pi^{st}_m(S^0)= \pi^{-m}(S^{0})$. 
 The $m=4$ case
is more subtle, however, and involves the $S^1$-equivariant machinery
in an essential manner. 

With Proposition \ref{woosh} and Corollary \ref{smoosh},
Theorem \ref{dimsum} implies 

\begin{prop} \label{squish} 
Let $X_j$, $j= 1, \ldots, 4$,   be smooth compact almost-complex
$4$-manifolds with non-zero  mod-$2$ Seiberg-Witten 
	invariant.  Let $N$ be a smooth compact oriented
$4$-manifold with $b_+(N)=0$, and
let $E_1 , \ldots , E_k$ be a set of generators for $H^2(N, {\mathbb Z})/\mbox{torsion}$
relative to which the   intersection form is diagonal. 
If $b_{1}(X_{j})=0$, $b_{+}(X_{j})\equiv 3\bmod 4$,
         and $\sum_{j=1}^4b_{+}(X_{j})\equiv 4 \bmod 8$, 
then, for any $m = 1 , \ldots, 4$,  
$$\sum_{j=1}^m \pm c_1(X_j) + \sum_{i=1}^k \pm  E_i$$
is a monopole class of
$[\#_{j=1}^m X_j ]\# N$. Here the $\pm$ signs are 
arbitrary, and are independent of one another. 
\end{prop}

\begin{cor} \label{squash} 
Let $X_j$ and $N$ be as in Proposition \ref{squish}, and, for
some  $m = 1 , \ldots, 4$,  
let $g$ be a Riemannian metric on 
$M=[\#_{j=1}^m X_j ]\# N$. 
Then, with respect to the Hodge decomposition 
$$H^2(M, \RR ) = {\mathcal H}^+_g \oplus {\mathcal H}^-_g , $$
 there is a monopole class 
$a\in H^2(M, \RR)$ whose self-dual part 
$a^+$ satisfies
$$(a^+)^2 \geq \sum_{j=1}^m c_1^2(X_j) .$$
\end{cor}

\begin{proof}
Let $E_1 , \ldots , E_k$ be a set of generators for $H^2(N, {\mathbb Z})/\mbox{torsion}$
relative to which the   intersection form is diagonal; such a set of generators always exists
\cite{donaldson,baufu}. 
Let $\alpha = \sum_{j=1}^m c_1(X_j)$, and choose new generators $\hat{E}_i= \pm E_i$
for $H^2(N, \ZZ)/\mbox{torsion}$ 
such that 
$$\alpha^+\cdot \hat{E}_i \geq 0 .$$
Then $a = \alpha + \sum_{i=1}^k\hat{E}_i$
is a monopole class on $M$, and 
\begin{eqnarray*}
(a^+)^2&=& \left[\alpha^+ + \sum_{i=1}^k \hat{E}_i^+\right]^2
\\&=&
(\alpha^+)^2 + 2 \alpha^+ \cdot \sum_{i=1}^k \hat{E}_i^+ + \left[\left(\sum_{i=1}^k \hat{E}_i\right)^+ \right]^2
\\  &=& 
(\alpha^+)^2 + 2  \sum_{i=1}^k ( \alpha^+\cdot \hat{E}_i) + \left[\left(\sum_{i=1}^k \hat{E}_i\right)^+ \right]^2
\\& \geq & (\alpha^+)^2
 \geq  \alpha^2
= \sum_{j=1}^m c_1^2(X_j), 
\end{eqnarray*}
as promised. 
\end{proof}

\section{Scalar Curvature Problems} 
\label{yada}

 Recall that  a  {\em conformal class} on  
 a smooth compact manifold  $M$
is by definition the set of 
smooth  Riemannian metrics on $M$ 
$$\gamma = [g]=\{ ug ~|~u: M\to {\mathbb R}^+\}$$ 
which are point-wise proportional to some Riemannian metric $g$.
As we saw in \S \ref{intro}, each  conformal class $\gamma$ has an  
associated number 
\begin{equation}\label{yconst}
{\mathcal Y}_{\gamma }
 = \inf_{{g}\in \gamma } \frac{\int_M 
s_{g}~d\mu_{g}}{\left(\int_M 
d\mu_{g}\right)^{\frac{n-2}{n}}},
\end{equation}
 called the
{\em Yamabe constant} of $\gamma$.
 For any metric $g$,  ${\mathcal Y}_{[g]}$
has the same sign as  the lowest eigenvalue of the {\em Yamabe Laplacian}
$\Delta + \frac{(n-2)}{4(n-1)} s_g$;  if $s_g$ does not change sign, moreover, 
this sign agrees with that of $s_g$. 
 A deep theorem
\cite{aubyam,lp,rick} 
of Yamabe, Trudinger, Aubin, and Schoen 
asserts that any conformal class 
$\gamma$ contains a 
metric  
which actually minimizes the relevant functional.  
Such a metric is called a {\em Yamabe minimizer}.
Any Yamabe minimizer has constant scalar curvature.
Conversely, any metric $g$ 
with $s_g=\mbox{const} \leq 0$ is 
 a Yamabe minimizer.

Given a smooth compact $n$-manifold $M$,
one 
 defines  the 
  {\em Yamabe invariant}  \cite{okob,sch,lno,lky}
 of $M$ by 
$${\mathcal Y}(M) = \sup_{\gamma\in {\mathcal C}(M)} {\mathcal Y}_{\gamma } = 
\sup_{\gamma}\inf_{{g}\in \gamma } \frac{\int_M
s_{g}~d\mu_{g}}{\left(\int_M 
d\mu_{g}\right)^{\frac{n-2}{n}}},$$
where ${\mathcal C}(M)$ denotes the space of all conformal classes
of metrics on $M$.  
As already pointed out in \S \ref{intro}, this is 
a real-valued diffeomorphism invariant of $M$. 
In light of the above comments, it  is immediate 
  that ${\mathcal Y}(M) > 0$ iff $M$ admits a 
metric of positive scalar curvature. On the other hand, 
if $M$ does not admit metrics of positive scalar curvature, the invariant
is simply the   supremum of 
the scalar curvatures of unit-volume 
constant-scalar-curvature metrics on $M$, 
since any constant-scalar-curvature metric
of non-positive scalar curvature is automatically a 
Yamabe minimizer. 

The minimax definition of ${\mathcal Y}(M)$ is technically rather  unwieldy.   Fortunately,
by an observation  of  Besson-Courtois-Gallot
\cite{bcg1}, 
${\mathcal Y}(M)$ and the conceptually simpler  invariant 
  $${\mathcal I}_s(M)= \inf_g \int_M|s_g|^{n/2}d\mu_g$$ 
  determine each other
  whenever   $M$ does not admit metrics of positive scalar curvature  \cite{and,lky}.

\begin{prop}\label{func} 
Let $M$ be a smooth compact $n$-manifold, $n \geq 3$. 
Then 
$$
{\mathcal I}_s (M) = \left\{
\begin{array}{ll} 0&\mbox{if } {\mathcal Y}(M)\geq 0\\
|{\mathcal Y}(M)|^{n/2} &\mbox{if } {\mathcal Y}(M)\leq 0.
\end{array}\right.
$$
\end{prop}
The proof depends on two main observations. First of all, 
there are
always conformal classes with negative Yamabe constant when 
$n\geq 3$;  because the  lowest eigenvalue of the Yamabe Laplacian
 depends continuously
on $g$,  $Y(M) > 0$ would therefore  imply the existence of a scalar-flat
metric, and thus force 
${\mathcal I}_s (M)$ to vanish. 
Secondly, 
detailed calculation reveals that the $L^{n/2}$ norm of $s$ 
is minimized in each conformal class 
of negative Yamabe constant by the metrics of constant scalar curvature.

This same argument also shows that  ${\mathcal I}_s$ 
 can be re-expressed \cite{bcg1,lky} as
$${\mathcal I}_s(M)  = \inf_g \int_M |s_{-g}|^{n/2}d\mu_g , $$
where $s_{-g}(x)= \min (s_{g}(x) , 0)$. 
Thus, in essence, we are allowed to 
neglect regions of positive scalar curvature when computing 
${\mathcal I}_s$. This observation has  some  important consequences. For
example, 
 if 
$(X,g_X)$ and $(Y,g_Y)$ are any compact Riemannian $n$-manifolds, 
then, for any $\varepsilon > 0$, the Gromov-Lawson method
\cite{gvln} of joining $X$ to $Y$ by a long, thin neck with $s >  0$ 
yields  a  Riemannian metric $g_{X\# Y}$ on 
the connected sum $X\# Y$ such that
$$
 \int_{X\# Y} |s_{-g_{X\# Y}}|^{n/2}d\mu_{g_{X\# Y}} < \int_X |s_{g_X}|^{n/2}d\mu_{g_X} + \int _Y|s_{g_Y}|^{n/2}d\mu_{g_Y} + \varepsilon . 
$$
One therefore has  the following general  inequality \cite{okob}:

\begin{prop} \label{cobble} 
Let $X$ and $Y$ be smooth compact $n$-manifolds, $n\geq 3$. Then 
$${\mathcal I}_s (X\# Y ) \leq {\mathcal I}_s (X ) + {\mathcal I}_s (Y).$$
\end{prop}

Now  the invariant ${\mathcal I}_s$ of 
a smooth compact $4$-manifold $M$   is just 
$${\mathcal I}_s(M) = \inf_g \int_M s^2_g d\mu_g , $$
so    Proposition  \ref{west} and Corollary \ref{squash} 
give us the following lower bound:\begin{prop} \label{bound}
Let $X_j$, $j= 1, \ldots, 4$,   be smooth compact almost-complex
$4$-manifolds with non-zero  mod-$2$ Seiberg-Witten 
	invariant.  Let $N$ be a smooth compact oriented
$4$-manifold with $b_+(N)=0$.  
If $b_{1}(X_{j})=0$, $b_{+}(X_{j})\equiv 3\bmod 4$,
         and $\sum_{j=1}^4b_{+}(X_{j})\equiv 4 \bmod 8$, 
then, for any $m = 1 , \ldots, 4$,  
$${\mathcal I}_s([\#_{j=1}^m X_j ]\# N)\geq 32\pi^2 \sum_{j=1}^m c_1^2(X_j).$$
\end{prop} 

Combining Propositions   \ref{cobble} and \ref{bound},  we thus obtain
\setcounter{main}{1}
\begin{main}
Let $X_j$, $j= 1, \ldots, 4$,   be minimal  compact complex surfaces with 
$b_{1}(X_{j})=0$, $b_{+}(X_{j})\equiv 3\bmod 4$,
         and $\sum_{j=1}^4b_{+}(X_{j})\equiv 4 \bmod 8$. 
  Let $N$ be a smooth compact oriented
$4$-manifold with $b_+(N)=0$ which admits a metric of 
non-negative  scalar curvature.  
Then 
$${\mathcal I}_s([\#_{j=1}^m X_j ]\# N)= 32\pi^2 \sum_{j=1}^m c_1^2(X_j).$$
\end{main} 
\begin{proof}
For any minimal compact complex surface $X$ with $b_+ > 1$, 
one has \cite{lno,lky} 
$${\mathcal I}_s (X) = 32\pi^2 c_1^2(X). $$
On the other hand, since $N$ admits a metric of non-negative scalar curvature,
$${\mathcal I}_s (N) =0.$$
Propositions \ref{cobble} and  \ref{bound}  thus tells us that
\begin{eqnarray*}
{\mathcal I}_s([\#_{j=1}^m X_j ]\# N)&\leq& {\mathcal I}_s (N)  + \sum_{j=1}^m {\mathcal I}_s (X_j)\\
&=& 32\pi^2 \sum_{j=1}^m c_1^2 (X_j) \\
&\leq& {\mathcal I}_s([\#_{j=1}^m X_j ]\# N), 
\end{eqnarray*}
and the promised equality  follows.
\end{proof} 

But if $M$ is a  $4$-manifold   
with a non-trivial Bauer-Furuta invariant, $M$ does not admit
metrics of positive scalar curvature, and hence has 
${\mathcal Y}(M)\leq 0$. 
Theorem A and Proposition \ref{func} therefore
together imply 

\setcounter{main}{0}
\begin{main}
Let $X_j$, $j= 1, \ldots, 4$,   be minimal  compact complex surfaces with 
$b_{1}(X_{j})=0$, $b_{+}(X_{j})\equiv 3\bmod 4$,
         and $\sum_{j=1}^4b_{+}(X_{j})\equiv 4 \bmod 8$. 
  Let $N$ be a smooth compact oriented
$4$-manifold with $b_+(N)=0$ 
and ${\mathcal Y}(N) \geq 0$. 
Then 
$${\mathcal Y}([\#_{j=1}^m X_j ]\# N)= -4\pi\sqrt{2 \sum_{j=1}^m c_1^2(X_j)}.$$
\end{main} 
 
\section{Ricci Curvature}

In dimension $n=4$, the Ricci-curvature-based invariant 
defined by equation (\ref{tight})  becomes
$$
{\mathcal I}_r(M)= \inf_g \int_M |r_g|^2d\mu_g
= \inf_g \int_M \Big(
\frac{s_g^2}{4}+ |\stackrel{\circ}{r}_g|^2\Big) d\mu_g , 
$$
where $\stackrel{\circ}{r}$ denotes the trace-free piece of the 
Ricci curvature. Now notice that  the  Gauss-Bonnet-type formula \cite{bes,hit,tho}
\begin{equation}\label{gb}
(2\chi + 3\tau) (M) =
\frac{1}{4\pi^2}\int_M \left(\frac{s_g^2}{24}+ 2|W_+|_g^2 - \frac{|\stackrel{\circ}{r}_g|^2}{2}\right) d\mu_g
\end{equation}
therefore tells us that 
\begin{equation}\label{gbric}
\int_M |r_g|^2d\mu_g = - 8\pi^2 (2\chi + 3\tau ) (M) +
 \int_M \left(\frac{s_g^2}{3}+ 4 |W_+|_g^2\right)d\mu_g
 \end{equation}
for 
any Riemannian metric $g$ on a compact oriented  $4$-manifold $M$.
On  the other hand, the Cauchy-Schwarz and triangle inequalities 
imply \cite{lric} 
that 
$$\int_M \left(\frac{s_g^2}{3}+ 4 |W_+|_g^2\right)d\mu_g
\geq\frac{2}{9} \int_M (s-\sqrt{6}|W_+|)|_g^2 d\mu_g . 
$$
If $a$ is a monopole class on $M$, Proposition \ref{est} therefore implies \cite{lric}
that 
\begin{equation}\label{ricci}
\int_M |r|_g^2 d\mu_g \geq 8\pi^2 \left[ 2(a^+)^2 - (2\chi + 3\tau ) (M)\right] .
\end{equation} 
We therefore have  
\begin{prop}
Let $X_j$, $j= 1, \ldots, 4$,   be smooth compact almost-complex
$4$-manifolds with non-zero  mod-$2$ Seiberg-Witten 
	invariant.  Suppose, moreover, 
that $b_{1}(X_{j})=0$, $b_{+}(X_{j})\equiv 3\bmod 4$,
         and $\sum_{j=1}^4b_{+}(X_{j})\equiv 4 \bmod 8$. 
 Let $N$ be a smooth compact oriented
$4$-manifold with $b_+(N)=0$. 
Then, for
  $m = 1 , \ldots, 4$,  
$${\mathcal I}_r (M) 
\geq    8\pi^2 \left[ 4m- (2\chi + 3\tau ) (N)+ 
\sum_{j=1}^m c_1^2 (X_j)\right] .$$
\end{prop}
\begin{proof}
By Corollary \ref{squash}, there is a 
monopole class $a$ on $M= [\#_{j=1}^m X_j ]\# N$
with $(a^+)^2 \geq  \sum_{j=1}^m c_1^2 (X_j)$. 
On the other hand,  
$$(2\chi + 3\tau ) (M)
=  
(2\chi + 3\tau ) (N)
-4m + 
\sum_{j=1}^m c_1^2 (X_j) , $$
so that  (\ref{ricci}) implies  that 
$$\int_M |r|_g^2 d\mu_g \geq 8\pi^2 \left[ 4m- (2\chi + 3\tau ) (N)+ 
\sum_{j=1}^m c_1^2 (X_j)\right] $$
for any metric $g$ on $M$. 
Taking the infimum over $g$ of the left-hand side now yields the 
desired inequality.
\end{proof} 

Amazingly, this estimate is sharp in many cases:

\setcounter{main}{2}
\begin{main} 
Let $X_j$, $j=1, \ldots, 4$ be minimal complex surfaces satisfying
(\ref{tom}--\ref{harry}).  Let $N$ be a smooth compact oriented
$4$-manifold which admits an anti-self-dual metric of positive scalar curvature. 
Then, for each choice of $m=1, \ldots, 4$,  
$$ 
{\mathcal I}_r ([\#_{j=1}^mX_j]\#N)
=8\pi^2\left[4m - (2\chi + 3\tau ) (N) + \sum_{j=1}^m c_1^2(X_j)\right]. 
$$
\end{main}
\begin{proof}
Let us first observe \cite{bourg,lsd} 
 that any anti-self-dual $4$-manifold $(N,h)$ of positive scalar curvature
satisfies $b_+(N)=0$.  We thus have 
$$\int_M |r|_g^2 d\mu_g \geq 8\pi^2 \left[ 4m- (2\chi + 3\tau ) (N)+ 
\sum_{j=1}^m c_1^2 (X_j)\right] $$
by the previous result. Of course, by (\ref{gbric}), 
this equivalent to the statement that 
$$ \int_M \left(\frac{s_g^2}{3}+ 4 |W_+|_g^2\right)d\mu_g
\geq 16\pi^2 \sum_{j=1}^m c_1^2 (X_j),$$
and what we need to show is that, for any $\varepsilon > 0$, there is a
metric $g_{\varepsilon}$  on $M$ with 
\begin{equation}\label{close} 
 \int_M \left(\frac{s_{g_{\varepsilon}}^2}{3}+ 4 |W_+|_{g_{\varepsilon}}^2\right)d\mu_{g_{\varepsilon}}
< \varepsilon  +  16\pi^2 \sum_{j=1}^m c_1^2 (X_j).
\end{equation}

Such a metric can be constructed as follows: 
choose $m$ points $p_1, \ldots , p_m \in N$, and let 
$u$ be the conformal Green's function of 
$\{ p_1, \ldots , p_m\}$:
$$(\Delta + \frac{s}{6}) u = \sum_{j=1}^m \delta_{p_j} .$$
For  any positive constant $c$, 
the metric $cu^2 h$ on $N-\{ p_1, \ldots  , p_m\}$
is then scalar-flat, anti-self-dual, and asymptotically
flat, with $m$ ends.  Equip the pluricanonical model $\check{X}_j$
of $X_j$ with a K\"ahler-Einstein orbifold
metric \cite{aubin,rkob,yau}, choose a point $q_j$ of $\check{X}_j$, remove
a small metric ball around $q_j$, and glue the boundary of
this ball to the $ j^{th}$  boundary component of $N$ minus $m$ balls 
centered at $p_1, \ldots , p_m$. As the metrics do not quite match, 
one must use a partition of unity to smooth things together, but
this can be done in a manner  such that the $L^2$ norms  of the
curvature tensor on the transition 
annuli are as small as desired \cite{lno,lric}. In the same way, 
replace a neighborhood of each orbifold point of $\check{X}_j$ by a gravitational
instanton \cite{kron}. 
Now the K\"ahler-Einstein orbifold metric on $\check{X}_j$ has
$$\int_{\check{X}_j} \left( \frac{s^2}{3}+4|W_+|^2\right) d\mu 
= \int_{\check{X}_j} \frac{s^2}{2} d\mu  = 16\pi^2 c_1^2(X_j),$$
whereas 
$N-\{ p_1, \ldots  , p_m\}$ and the gravitational instantons
have  
$s\equiv 0$ and $W_+\equiv 0$.
The metric $g_{\varepsilon}$ can thus be chosen  
to satisfy (\ref{close}), and the
 result therefore follows. 
\end{proof}

In particular, since 
$k\cpb \# \ell \sxs$ admits anti-self-dual metrics
of positive scalar curvature \cite{jongsu,mcp2}, we obtain

\begin{cor}
Let $X_j$, $j=1, \ldots, 4$ be minimal complex surfaces satisfying
(\ref{tom}--\ref{harry}). Then, for any integers $k, \ell \geq 0$,
and any $m=1, \ldots, 4$, 
$$ 
{\mathcal I}_r \left([\#_{j=1}^mX_j]\#k\cpb \# \ell \sxs \right)
=8\pi^2\left[ k+ 4(\ell + m - 1)   + \sum_{j=1}^m c_1^2(X_j)\right]. 
$$
\end{cor}

\section{Obstructions to Einstein Metrics}\label{nein}

The main point of this article has been that
the Bauer-Furuta invariant can be used  to deduce certain
sharp curvature-integral estimates for many  interesting $4$-manifolds. 
However,  the
same techniques allow one to derive other 
 curvature estimates which,
while presumably not sharp, nonetheless have interesting consequences. 
In particular, 
Proposition \ref{est} can also be used to construct many new examples of 
smooth compact $4$-manifolds which do not admit  Einstein metrics. 

Indeed, the Cauchy-Schwarz and triangle inequalities 
 imply \cite{lric} 
that 
$$ \int \left(\frac{s^2}{24}+ 2 |W_+|^2\right)d\mu \geq 
\frac{1}{27} 
\int (s-\sqrt{6}|W_+|)^2 d\mu.$$
If $a$ is a monopole class on $M$, Proposition \ref{est} therefore tells us \cite{lric}
that 
\begin{equation}\label{einstein}
\frac{1}{4\pi^2}
\int \left(\frac{s^2}{24}+ 2 |W_+|^2\right)d\mu \geq \frac{2}{3}(a^+)^2.
\end{equation} 
By Corollary \ref{squash}, 
we therefore have

\begin{main} 
Let $X_j$, $j=1, \ldots, 4$ be  smooth, compact 
almost-complex $4$-manifolds for which the 
 mod-2 Seiberg-Witten invariant is non-zero, and 
suppose, moreover, that  
(\ref{tom}--\ref{harry}) hold.  Let $N$ be any oriented $4$-manifold
with $b_+=0$.  Then, for any $m=2, 3$ or    $4$,  the smooth $4$-manifold 
$[ \#_{j=1}^mX_j]\#N$ does not admit 
Einstein metrics if 
$$
 4m -   (2\chi + 3\tau )(N)    \geq \frac{1}{3}\sum_{j=1}^m c_1^2(X_j) .
$$
\end{main}
\begin{proof}
By Corollary \ref{squash}, 
for any choice of Riemannian metric $g$ on 
$M= [\#_{j=1}^m X_j ]\# N$, there is a 
monopole class $a$ on $M$ 
with $(a^+)^2 \geq  \sum_{j=1}^m c_1^2 (X_j)$. 
Thus any metric $g$ on $M$ satisfies 
$$\frac{1}{4\pi^2}
\int \left(\frac{s^2}{24}+ 2 |W_+|^2\right)d\mu \geq \frac{2}{3}  \sum_{j=1}^m c_1^2 (X_j)$$
by (\ref{einstein}),
and the inequality is necessarily strict if  
$m > 1$, since 
$M$ does not then admit a symplectic structure.
If $g$ is Einstein, however, the left-hand side
equals $(2\chi + 3\tau ) (M)$ by (\ref{gb}). Since
$$	(2\chi + 3\tau ) (M) =-4m+  (2\chi + 3\tau ) (N) +  \sum_{j=1}^m c_1^2 (X_j),$$
this means that the existence of an Einstein metric $g$ on $M$  implies,
for $m=2,3$ or $4$, that 
$$-4m+  (2\chi + 3\tau ) (N) +  \sum_{j=1}^m c_1^2 (X_j) > \frac{2}{3}  \sum_{j=1}^m c_1^2 (X_j),$$
or in other words that
$$
\frac{1}{3}  \sum_{j=1}^m c_1^2 (X_j) > 4m - (2\chi + 3\tau ) (N).
$$
The desired result therefore follows by contraposition. 
\end{proof}

	In particular,  using a celebrated result of Taubes \cite{taubes},
this immediately gives us  many new spin examples with free fundamental groups:

	\begin{cor}\label{non}
	Let $X$ be a simply connected  symplectic $4$-dimensional spin manifold
	with  $b_{+}\equiv 3 \bmod 8$. 
	Then 
  the spin manifold 
	 $X\# n K3 \# \ell \sxs$ does not
	admit Einstein metrics if $n= 1,2$ or $3$ and $\ell + n \geq c_1^2(X)/12$.
		\end{cor}

We remark that the work of Gompf \cite{gompf}
provides  an enormous catalog of choices for $X$. 
Moreover, by replacing a $K3$ with one of Kodaira's homotopy $K3$'s, 
one can find infinitely many smooth structures on these manifolds for which 
the obstruction also applies; cf. \cite{il1}. On the other hand, if we take $X$
to have negative signature, the Hitchin-Thorpe inequality \cite{bes,hit,tho} would only 
provide an obstruction to the existence of Einstein metrics on 
these spaces for $\ell + n \geq c_1^2(X)/4$.

\bigskip
\noindent
{\bf Acknowledgment.} 
The authors would like to thank Stefan Bauer and Mikio Furuta 
for many helpful remarks regarding  the  Bauer-Furuta invariant.

\vfill

{\footnotesize 
\noindent
{Masashi Ishida,
{Department of Mathematics,  
Sophia University, 
 Tokyo, Japan 102-8554}\\
{\sc e-mail}: ishida@mm.sophia.ac.jp}
 \\ {
Claude LeBrun,
Department of Mathematics, SUNY, 
Stony Brook, NY 11794-3651, USA\hfill 
{\sc e-mail}: claude@math.sunysb.edu}}
 
 \end{document}